\documentclass[11pt,a4paper,reqno]{amsart}

\usepackage[includehead,includefoot,margin=25mm]{geometry}

\usepackage{amsmath,amssymb}
\usepackage{mathrsfs}
\usepackage{enumerate}

\theoremstyle{plain}
\newtheorem{theorem}{Theorem}
\newtheorem{lemma}[theorem]{Lemma}

\theoremstyle{definition}
\newtheorem{definition}{Definition}
\newtheorem{example}{Example}

\theoremstyle{remark}
\newtheorem{remark}{Remark}

\newcommand{\wD}{\widehat{\mathscr{D}}}
\newcommand{\D}{{\mathscr{D}}}

\newcommand{\CC}{{\mathbb C}}
\newcommand{\QQ}{{\mathbb Q}}
\DeclareMathOperator{\Aut}{Aut}
\DeclareMathOperator{\Spec}{Spec}
\DeclareMathOperator{\divi}{div}
\newcommand{\C}{\mathbb{C}}
\newcommand{\Z}{\mathbb{Z}}
\newcommand{\Q}{\mathbb{Q}}
\newcommand{\A}{\mathbb{A}}

\newcommand{\GA}[0]{\ensuremath{\mathbb{G}_{\mathrm{a}}}}
\newcommand{\GM}[0]{\ensuremath{\mathbb{G}_{\mathrm{m}}}}

\newcommand{\spec}[0]{\ensuremath{\operatorname{Spec}}}
\newcommand{\OO}[0]{\ensuremath{\mathcal{O}}}

\begin{document}

\title[Characterization of Danielewski surfaces]{On the characterization of Danielewski surfaces
  by their automorphism groups}
  
\author{Alvaro Liendo}
\author{Andriy Regeta}  \author{Christian Urech}
  \thanks{Partially supported by Fondecyt projects 1160864 and 1200502.}
\address{Instituto de Matem\'atica y F\'\i sica \\ Universidad de Talca \\ Casilla 721, Talca, Chile}
\email{aliendo@inst-mat.utalca.cl}

\address{Institut f\"{u}r Mathematik \\ Friedrich-Schiller-Universit\"{a}t Jena, \\
Jena 07737, Germany}
\email{andriyregeta@gmail.com}
\thanks{Supported by the SNF, project number P2BSP2\_175008.}
\address{EPFL SB MATH \\ 
	Station 8 \\
	1015 Lausanne, Switzerland}
\email{christian.urech@gmail.com}

\maketitle

\begin{abstract}
In this note we show that if the automorphism group of a normal affine surface $S$ is isomorphic to the automorphism group of a Danielewski surface, then $S$ is isomorphic to the normalization of a Danielewski surface.
\end{abstract}

\section*{Introduction}
  Throughout this note we work over the field of complex numbers $\mathbb{C}$ and algebraic varieties are always considered to be affine and irreducible. One of our main results in \cite{liendo2018characterization} is the proof that affine toric surfaces are uniquely determined by their automorphism groups in the category of normal affine surfaces. In this note we apply similar techniques to investigate in as far this result can be extended to other classes of affine surfaces with a large automorphism group.
  
  A well studied class of affine surfaces are {\it Danielewski surfaces}, i.e., surfaces of the form $\D_p^n = \{ x^ny = p(z) \} \subset \mathbb{A}^3$ for some polynomial $p\in \C[z]$. We denote by $\wD_p^n$ the normalization of $\D_p^n$. Recall that $\D_p^1$ is always normal.  
  These surfaces were introduced by Danielewski in order to construct a counterexample to the generalized Zariski cancellation problem (\cite{danielewski1989cancellation}). Since then, numerous papers have been published on the subject, in particular with regards to the rich structure of their automorphism groups.  The automorphism groups of two generic smooth Danielewski surfaces $\D_p^1$ and $\D_q^1$ are isomorphic, where $\D_p^1$ is generic if $\deg p(z)\geq 3$ and no affine automorphism permutes the roots of $p(z)$ in $\CC$. This follows from \cite[Theorem and Remark~(3) on page 256]{MR1045377}, and more precisely from \cite[Theorem~2.7]{MR3495426}. Indeed, in this last reference, it is proven that for a generic Danielewski surface $\D_p^1$, we have  $\Aut(\D_p^1)\simeq (\CC[x]\ast\CC[y])\rtimes(\CC^*\rtimes \Z/2\Z)$ and the semidirect product structure does not depend on $p(z)$.
 A similar result holds for $n>1$ where we have $\Aut(\D_p^n)\simeq (\CC[x])\rtimes\CC^*$ for every generic polynomial $p$, see \cite{makar2001group}.  
  This yields that the automorphism group of a
  Danielewski surface does not determine the surface in general. However, we prove the following result:
  
  \begin{theorem}\label{Danielewski}
    Let $S$ be a normal affine surface and $\wD_p^n$ be the normalization of a Danielewski surface $\D_p^n$ for some $p\in\C[z]$ and $n\in \Z_{>0}$. If $\Aut(S)$ and $\Aut(\wD_p^n)$ are isomorphic as groups, then $S$ is isomorphic to the normalization $\wD_q^m$ of a Danielewski surface $\D_q^m$ for some polynomial $q \in \mathbb{C}[z]$ and some $m\in \Z_{>0}$. Moreover, if $n=1$ then $S$ is isomorphic to $\D_q^1$ for some polynomial $q \in \mathbb{C}[z]$.
  \end{theorem}
  
  Let $\GM$ and $\GA$ be the multiplicative and the additive group over $\CC$, respectively. All Danielewski surfaces $\D_p^n$, and hence their normalizations $\wD_p^n$, admit a $\GM$-action given in the ambient space $\mathbb{A}^3$ via $t\colon (x,y,z)\mapsto (tx,t^{-n}y,z)$ for $t\in\GM$. The main idea of the proof of Theorem~\ref{Danielewski} is contained in Lemma~\ref{gap-FZ}, which characterizes $\wD_p^n$ in terms of certain extensions of  this $\GM$-action by $\GA$.
 
  \begin{remark}
    It is proved in \cite{Leuenberger:aa} that for two polynomials $p$ and $q$ with simple roots, $\Aut(\D_p^1)$ is isomorphic to $\Aut(\D_q^1)$ as a so-called {\it ind-group} if and only if $\D_q^1$ is isomorphic to $\D_p^1$ as a variety, as opposed to the case of abstract group isomorphisms. The main reason for this comes from the additional rigidity of Lie algebras. Indeed, the Lie algebra of $\Aut(\D_p^1)$ is isomorphic to the Lie algebra of $\Aut(\D^1_q)$ if and only if $\D^1_p$ is isomorphic to $\D^1_q$ and ind-group isomorphisms induce isomorphisms of the corresponding Lie algebras. Together with Theorem \ref{Danielewski} this gives us that a surface isomorphic to $\D_p^1$ is determined by its automorphism group seen as an ind-group in the category of smooth affine surfaces.
  \end{remark}

  \subsection*{Acknowledgements} The authors would like to thank M.\,Zaidenberg for useful discussions and both anonymous referees for numerous comments that helped to improve the exposition of the paper. Part of this work was done during a stay of the three authors at IMPAN in Warsaw. We would like to thank IMPAN and the organizers of the Simons semester ``Varieties: Arithmetic and Transformations'' for the hospitality. The second author {would like to thank the Max-Plank Institute for Mathematics in Bonn for its kind support.  }
	
  \noindent This work was partially supported by the grant 346300 for IMPAN from the Simons Foundation and the matching 2015-2019 Polish MNiSW fund.  

 \section*{Root subgroups of non-toric
   $\GM$-surfaces} \label{sec:non-toric}

Let $S$ be a an affine surface and $G$ be an algebraic group. A regular faithful action of $G$ on $S$ induces an injective homomorphism from $G$ to $\Aut(S)$. We say that the image of $G$ in $\Aut(S)$ is an \emph{algebraic subgroup} of $\Aut(S)$. One can show that an algebraic subgroup admits a canonical structure of algebraic variety \cite[Theorem~0.3.1]{furter2018geometry}.

 A $\GM$-surface is a surface $S$ together with a given regular faithful $\GM$-action on $S$. Let $T\subset \Aut(S)$ be the acting torus, i.e., the image of $\GM$ in $\Aut(S)$. A {\it root subgroup} of $S$ with respect to $T$ is an algebraic subgroup $U\subset\Aut(S)$ isomorphic to $\GA$ that is normalized by $T$. Let $\lambda\colon \GA\to U$ be an isomorphism. There exists a character $\chi\colon T\to\GM$ not depending on the choice of $\lambda$ such that $t\circ\lambda(s)\circ t^{-1}=\lambda(\chi(t)s)$. This character is called the \emph{weight} of $U$.  Recall that the set of characters $\chi\colon T\rightarrow \GM$ forms a group $\mathfrak{X}(T)$ isomorphic to $\Z$ and such an isomorphism is uniquely determined up to sign. In \cite{MR2020670} a classification of normal affine $\GM$-surfaces was given, followed by a classification of their root subgroups in \cite{MR2196000}. We recall here the main features of the classification that we need in this paper. This is a short version of our account of the subject in \cite[Section~4]{liendo2018characterization}.

 \begin{definition}
    Surfaces endowed with a $\GM$-action are classified in three \emph{dynamical types} \cite{MR0460342}: a $\GM$-surface is \emph{elliptic} if the $\GM$-action has an attractive fixed point, \emph{parabolic} if the $\GM$-action has infinitely many fixed points and \emph{hyperbolic} if the $\GM$-action has at most finitely many fixed points none of which is attractive. 
 \end{definition}

A $\GM$-action $\alpha\colon\GM\times S\rightarrow S$ on an affine surface $S$ induces a $\mathfrak{X}(T)$-grading on the algebra of regular functions. Under the isomorphism $\mathfrak{X}(T)\simeq \mathbb{Z}$, it is customary to denote this as a $\Z$-grading of the algebra of regular functions given by
 \[\OO(S)=\bigoplus_{i\in \Z}A_i, \quad\mbox{where}\quad A_i=\left\{f\in\OO(S)\mid \alpha^*(f)=t^i\cdot f \right\}\,.\]

 The elements in $A_i$ are called \emph{semi-invariants of weight $i\in \Z$}. A $\GM$-surface is hyperbolic if and only if there exist non-trivial semi-invariants whose weights have different sign. In the hyperbolic case, generic orbit closures are isomorphic to $\A^1_*$. If the surface is not hyperbolic, all semi-invariants that are not invariant have the same sign. In this case,  the normalizations of the generic orbit closures are isomorphic to $\A^1$. The elliptic case corresponds to the case where the only invariant functions are the constants. The parabolic case corresponds to case where the ring of invariant functions has transcendence degree 1 over $\CC$ and therefore there is a curve of  points fixed by $\GM$ in the surface.

 \medskip

 In algebraic terms, root subgroups are in one to one correspondence with \emph{homogeneous locally nilpotent derivations} of the $\Z$-graded algebra $\OO(S)$. A { homogeneous} locally nilpotent derivation is a $\C$-linear map $\delta\colon\OO(S)\rightarrow\OO(S)$ that sends semi-invariants to semi-invariants, satisfies the Leibniz rule $\delta(fg)=f\delta(g)+g\delta(f)$ for all $f,g\in \OO(S)$ and for every $f\in \OO(S)$ there exists $n\in \Z_{>0}$ such that $\delta^n(f)=0$, where $\delta^n$ denotes the composition of $\delta$ with itself $n$-times. In particular, the Leibniz rule implies that for every homogeneous locally nilpotent derivation $\delta$ there exists an integer $\ell$ such that $\delta(A_i)$ is contained in $A_{i+\ell}$ for any $i \in \mathbb{Z}$. We call $\ell$ the \emph{degree} of $\delta$.  Recall that under the isomorphism $\mathfrak{X}(T)\simeq \mathbb{Z}$, the degree of $\delta$ corresponds to a character of the acting torus $T$.
 See \cite[Section~4.1]{liendo2018characterization} for a more detailed description of root subgroups in terms of homogeneous locally nilpotent derivations. The next theorem summarizes the results from \cite{liendo2018characterization} as needed for this paper.

  \begin{theorem}\label{semidirect-product}
    Let $S$ and $S'$ be normal surfaces with $S$ non-toric. Assume that $\Aut(S)$ contains algebraic subgroups $T$ and $U$ isomorphic to $\GM$ and $\GA$, respectively. Let $\varphi\colon \Aut(S)\rightarrow\Aut (S')$ be a group isomorphism, then the following hold:
     \begin{enumerate}
     \item The image $\varphi(T)\subset\Aut(S')$ is an algebraic subgroup isomorphic to $\GM$. \label{torus}
     \item There exist root subgroups in $\Aut(S)$ and they are mapped to root subgroups preserving weights, up to a torus isomorphism not depending on the root subgroup. \label{roots}
     \item The surfaces $S$ and $S'$ are of the same dynamical type. \label{types}
     \end{enumerate}
  \end{theorem}

\begin{proof}
The statements (a) and (b) are proven in \cite[Theorem~6.5]{liendo2018characterization}. Statement (c) follows directly from \cite[Theorem~1.2]{liendo2018characterization}.

\end{proof}

  We will also need the following lemma proven in \cite{liendo2018characterization}.

  \begin{lemma}[{\cite[Lemma~4.16]{liendo2018characterization}}]
    \label{different-weight} %
    A non-toric $\GM$-surface $S$ admits root subgroups of different weights if and only if $S$ is hyperbolic. Furthermore, in this case all root subgroups have different weights.
  \end{lemma}

  The following theorem borrowed from \cite[Section 4.2]{MR2020670} is the main classification result for hyperbolic $\GM$-surfaces.
  \begin{theorem} %
  Every hyperbolic affine $\GM$-surface is equivariantly isomorphic to $S=\spec A$, where
  $$A=\bigoplus_{i<0}H^0(C,\OO(\lfloor-iD_-\rfloor))\oplus
  \bigoplus_{i\geq0} H^0(C,\OO(\lfloor iD_+\rfloor))\,,$$ 
  where $C$ is the algebraic quotient of $X$ by $\GM$ and $D_+,D_-$ are two $\Q$-divisors on $C$ satisfying $D_++D_-\leq 0$.  Moreover, $S$ is uniquely determined by $C$ and the couple $(D_+,D_-)$ up to linear equivalence. In other words, the couples of divisors $(D_+,D_-)$ and $(D'_+,D'_-)$ on $C$ give rise to equivariantly isomorphic $\GM$-surfaces if and only if $D_+=D'_++\divi(h)$ and $D_-=D'_--\divi(h)$, for some rational function $h$ on the curve $C$.
  \end{theorem}

  \begin{example} \label{example} %
    In \cite[Example~4.10]{MR2020670} it is proven that the normalization $\wD_p^n$ of the Danielewski surface $\D_p^n$ is given by the data $D_+=0$ and $D_-=-\tfrac{1}{n}\operatorname{div}(p)$ on $C=\A^1$. Remark that any $\QQ$-divisor $D$ in $\A^1$ with negative coefficients gives rise to a normalization of a Danielewski surface by taking $D_+=0$ and $D_-=D$.
\end{example}

\begin{lemma} \label{gap-FZ-1} %
  Let $S$ be a non-toric $\GM$-surface that is given by the couple of divisors $(D_+,D_-)$ in $C$. If there are two root subgroups with non-negative weights with respect to $T$ in $\Aut(S)$ whose weights differ by one then $C=\A^1$ and $D_+$ is integral.
\end{lemma}

\begin{proof}
  Assume there are two root subgroups with respect to $T$ in $\Aut(S)$ whose weights differ by one. Since there exist root subgroups of different weights, by Lemma~\ref{different-weight}, we have that $S$ is hyperbolic. In the language of \cite{MR2020670}, $S$ is described by the couple of $\Q$-divisors $D_+$ and $D_-$ on a smooth affine curve $C$. By \cite[Theorem~3.22]{MR2196000}, we have $C\simeq \A^1$, and up to linear equivalence, we can assume $D_+=-\tfrac{e'}{d}\cdot[0]$ and the weight $e$ of a root subgroup must satisfy $ee'=1\mod d$, where $0\leq e'<d$. But $S$ admits root subgroups whose weights differ by 1. This yields $d=1$ and so $e'=0$ or, equivalently, $D_+=0$ up to linear equivalence. This yields the lemma. 
\end{proof}

From Lemma~\ref{gap-FZ-1} and Example~\ref{example} we deduce the following criterion characterizing normalizations of Danielewski surfaces among normal $\GM$-surfaces.

\begin{lemma} \label{gap-FZ} %
  Let $S$ be a non-toric $\GM$-surface. Then there are two root subgroups with respect to $T$ in $\Aut(S)$ whose weights differ by one if and only if $S$ is the normalization of a Danielewski surface.
\end{lemma}
\begin{proof}
 Up to torus automorphism, we assume that both root subgroups have non negative weights. Now by Lemma~\ref{gap-FZ-1} we obtain that $D_+$ is integral and so up to linear equivalence, we can assume that $D_+=0$. 
 Now, by Example~\ref{example}, it follows that $S$ is isomorphic to the normalization of a Danielewski surface.
 
  On the other hand, the Danielewski surface $\D_p^n$ admits a root subgroup $U$ given by the homogeneous locally nilpotent derivation $\delta$ given by $\delta(x)=0$, $\delta(y)=p'(z)$ and $\delta(z)=x^n$ whose weight is $n$. Furthermore, since $x$ is $U$-invariant and $\GM$-semi-invariant of weight 1, we conclude that $x\delta$ is also a homogeneous locally nilpotent derivation and its corresponding root subgroup has weight $n+1$. This proves the lemma since connected algebraic group actions lift to the normalization by the universal property.
\end{proof}

\begin{remark} \label{exa-dp1} %

  For the proof of Theorem~\ref{Danielewski} we need to compute all the weights of root subgroups in $\Aut(\D_p^1)$. In general, the normalization $\wD_p^n$ of the Danielewski surface $\D_p^n$ is given by the data $D_+=0$ and $D_-=-\tfrac{1}{n}\operatorname{div}(p)$ on $\A^1$ as in Example~\ref{example}. Applying \cite[Theorem~3.22]{MR2196000} a routine computation yields that the non-negative weights of root subgroups in $\Aut(\wD_p^n)$ are exactly the integers greater than or equal to $\tfrac{n}{l}$, where $l$ is the smallest order of a root of $p(z)$. Hence, in the case of $\Aut(\D_p^1)$ we obtain that all positive integers appear as weights of root subgroups. 
  
  To compute the negative numbers that appear as weights of root subgroups, we reverse the grading taking the automorphism $t\mapsto t^{-1}$ of $\GM$. This accounts to exchanging $D_+$ and $D_-$.  Now a similar application of \cite[Theorem~3.22]{MR2196000} yields that all negative numbers appear as the weight of a root subgroup in $\Aut(\D_p^1)$.
\end{remark}

In the following lemma we show that the normalizations of $\D_p^n$ and $\D_{p^d}^{dn}$ coincide.

\begin{lemma} \label{normal-Danielewski} %
  Let $\D_p^n$ and $\D_q^m$ be two Danielewski surfaces and let $\wD_p^n$ and $\wD_q^m$ be their normalizations respectively, where $p(z),q(z)\in \C[z]$ and $n,m\in \Z_{>0}$. If $n=d\cdot m$ and $p(z)=q^d(z)$ for some $d\in\Z_{>0}$ then $\wD_p^n$ and $\wD_q^m$ are isomorphic.
\end{lemma}

\begin{proof}
  Let $A=\OO(\D_p^n)=\C[x,y,z]/(x^ny-p(z))$. The field of rational functions of $A$ is $\C(x,z)$. The element $f=\tfrac{q(z)}{x^m}$ belongs to the normalization $\widetilde{A}$ of $A$ since it satisfies the equation $y-f^d=0$. This yields an inclusion of algebras $A\subseteq B\subseteq \widetilde{A}$, where $B=\C[x,z,f]$. The lemma follows since $B\simeq \OO(\D_q^m)$.
\end{proof}

\begin{lemma} \label{lemma-dp1} %
  Assume $\wD_p^n$ has root subgroups of all weights different from zero. Then $\wD_p^n$ is isomorphic to $\D_q^1$ for some $q\in \C[z]$.
\end{lemma}

\begin{proof}
  The surface $\wD_p^n$ is given by the combinatorial data $D_+=0$ and $D_-=-\tfrac{1}{n}\operatorname{div}(p)$ in the algebraic quotient $\A^1=\Spec\C[z]$ of the $\GM$-action. By reversing the grading we exchange the roles of $D_+$ and $D_-$. Since there are two root subgroups with non negative weights for the reverse grading, by Lemma~\ref{gap-FZ-1} we obtain that $D_-$ is integral, which is equivalent to the fact that there exists $q(z)\in \C[z]$ such that $D_-+\divi(q)=0$. It now follows that $D_-=-\tfrac{1}{n}\divi(p)=-\operatorname{div}(q)$.  This is equivalent to $p(z)=q^n(z)$. Finally, since $p(z)$ is a regular function on $\mathbb{A}^1$ the same holds for $q(z)$ and so $q(z)\in\CC[z]$. By Lemma~\ref{normal-Danielewski} we conclude that $\wD_p^n$ is isomorphic to $\wD_q^1$, the normalization of $\D_q^1$. A straightforward computation shows that $\D_q^1$ has only isolated singularities and so $\D_q^1$ is already normal by \cite[Chapter~2, Proposition~8.23]{MR0463157}. This concludes the proof.
\end{proof}

\section*{Proof of  Theorem~\ref{Danielewski}}

If $\wD_p^n$ is toric, then the result follows directly from \cite[Theorem~1.3]{liendo2018characterization}. In the sequel, we assume that $\wD_p^n$ is not toric. Let $\varphi\colon \Aut(\wD_p^n) \rightarrow \Aut(S)$ be an isomorphism of groups and let $T\subset\Aut(\wD_p^n)$ be the acting torus coming from the $\GM$-surface structure on $\wD_p^n$. Since $\wD_p^n$ is a hyperbolic $\GM$-surface, by Theorem~\ref{semidirect-product}, $S$ is hyperbolic, $\varphi(T)$ is an algebraic $1$-dimensional torus and root subgroups are mapped to root subgroups with the same weight up to torus automorphism. By Lemma~\ref{gap-FZ} there are two root subgroups with respect to $T$ in $\Aut(\wD_p^n)$ whose weights differ by one. We conclude that there are two root subgroups with respect to $\varphi(T)$ in $\Aut(S)$ whose weights differ by one. Again by Lemma~\ref{gap-FZ} we conclude that $S$ is isomorphic to the normalization of a Danielewski surface $\D_q^m$.

To prove the last statement of the theorem, recall {first} that $\D_p^1$ is always normal.  If $\Aut(S)$ is isomorphic to $\Aut(\D_p^1)$ then $S$ is isomorphic to the normalization of a Danielewski surface. Since the isomorphism $\varphi\colon \Aut(\D_p^1)\rightarrow \Aut(S)$ preserves weights of root subgroups, by Remark~\ref{exa-dp1} we have that $\Aut(S)$ has root subgroups with all possible non-zero weights. It follows from Lemma~\ref{lemma-dp1} that $S$ is isomorphic to $\D_q^1$ for some $q(z)\in \C[z]$. \qed

\end{document}